# The Asymmetric Traveling Salesman Problem

## Howard Kleiman

**I. Introduction.** Starting with M(a), an $n \times n$ asymmetric cost matrix, Jonker and Volgenannt [1] transformed it into a $2n \times 2n$ symmetric cost matrix, M(s), where M(s) has unusual properties. One such property is that an optimal tour in M(s) yields an optimal tour in M(a). Modifying M(s), we apply the modified Floyd-Warshall algorithm given in [2] to M(s). Let $T$ be a tour that is an upper bound for an optimal tour in M(a). Due to the structure of M(s), we either can always obtain an optimal tour in M(s) that is derived from only one minimal positively-valued cycle in $\sigma_T^{-1} M^-$ whose value is less than $|T|$, (i.e., we don't have to link circuits), or else $T = T_{OPT}$. Thus, we can obtain an optimal tour in M(a) in at most polynomial running time. If the proof of theorem 1 in section II is correct, since the asymmetric traveling salesman problem is NP-hard, P would equal NP.

## II. A Theorem

**Theorem 1.**

Let

$$M(a) = \begin{pmatrix} a_{11} & \cdots & a_{1n} \\ \vdots & \ddots & \vdots \\ a_{n1} & \cdots & a_{nn} \end{pmatrix}$$

be an asymmetric cost matrix where $a_{ii} = \infty, i = 1, 2, \ldots, n$. Using a modified version of the symmetric cost matrix, M(s), obtained by Jonker and Volgenannt [1] as well as a result of Kleiman in [2] and the use of the modified F-W algorithm, we prove that we always can obtain an optimal solution to M(a) in polynomial time.

**Proof.** If M(a) contains a non-positive entry, let $m$ be the smallest value of all the entries in M(a). We then add $-m + 1$ to each entry of M(a). Thus, each of the entries in M(a) now has a positive value. Jonker and Volgenannt [1] gave a method for transforming an $n \times n$ asymmetric cost matrix into a $2n \times 2n$ symmetric cost matrix such that an optimal tour in the latter yields an optimal tour in the former. Let $M_\infty$ be an $n \times n$ matrix each of whose entries is $\infty$. Furthermore, we change each diagonal entry of M(a) into $0$ to obtain the matrix $M(a)_d$. Finally, we define $M(a)_d^T = \begin{pmatrix} a_{11} & \cdots & a_{n1} \\ \vdots & \ddots & \vdots \\ a_{1n} & \cdots & a_{nn} \end{pmatrix}$. We then define the $2n \times 2n$ symmetric matrix

$\begin{pmatrix} M_\infty & M(a)_d \\ M(a)_d^T & M_\infty \end{pmatrix}$ as M(s). We use any algorithm that yields an upper bound, say $T_{UPPERBOUND} = (t_1\ t_2\ \ldots\ t_n)$, for an optimal tour in M(a). In M(s), $a_{i+n, i} = a_{n+i, i} = 0, i = 1, 2, \ldots, n$. We now replace $T_{UPPERBOUND}$ by

$T = (t_1 \ t_{2+n} \ t_2 \ t_{3+n} \ t_3 \ ... \ t_n \ t_{1+n})$ in M(s). By construction, $|(t_{i+n} \ t_i)| = 0, \ i = 1, 2, \ ... \ , n$ in M(s). It follows that $\prod_{i=1}^{i=n}(t_{i+n} \ t_i) = \sigma_T$ where each 2-cycle $(t_{i+n} \ t_i)$ has a value of $0$. We can always use a product of 2-cycles (edges) to obtain $\sigma_T^{-1} M(s)^-$. As mentioned earlier, the J-V paper proves that an optimal tour of M(s) yields an optimal tour of M(a). An acceptable path in M(s) consists alternately of non-zero and zero arcs. We cannot link acceptable cycles of the kind found in M(s) since by linking by deleting two arcs of form $(t_{i+n} \ t_i)$ or $(t_i \ t_{i+n})$, we obtain a circuit containing two consecutive *non-zero-valued* directed edges. Using the modified F-W algorithm, each cycle from $a$ to $b$ obtained in M(s) has a value no greater than any other cycle from $a$ to $b$. Thus, - using the modified F-W algorithm -, there is only one way that we could obtain $T_{FWOPT}$ of M(s): one minimal positively-valued acceptable cycle containing $n$ arcs whose value is less than $|T|$ or – if one can't be found, $T_{UPPERBOUND} = T_{FWOPT} = T_{OPT}$. As proved in [2], this cycle always yields an *optimal* tour in M(s) that yields an optimal tour in M(a). We now show that the modified F-W algorithm when used for obtaining acceptable paths always obtains all acceptable paths in at most $O(n^4)$ running time. Since each such cycle is obtained using the modified F-W algorithm – together with an algorithm to insure that an acceptable path obtained stays acceptable which requires backtracking in only a smaller number of cases than otherwise – we can obtain a such a minimal positively-valued acceptable cycle containing $n$ points of value less than $|T|$ (if it exists) in polynomial time. In particular, the Floyd-Warshall algorithm has $O(n^3)$ running time. Thus, even backtracking in every case, would raise the running time to at most $O(n^4)$.

### III. The Construction of $\sigma_T^{-1}(M(s))^-$

$$\begin{pmatrix} a_{11} & \cdots & a_{1n} \\ \vdots & \ddots & \vdots \\ a_{n1} & \cdots & a_{nn} \end{pmatrix} \quad M(a)$$

$a_{ii} = \infty, \ i = 1, 2, \ ... \ , n$. In J-V M(a), $a_{ii} = -M'$ where $M'$ is the largest value of a non-diagonal entry in M(a). In $M(a)_d$, $a_{ii} = 0$. $M(a)_d$ is used in the construction of M(s). In order that an optimal tour of M(s) yields an optimal tour of M(a), all arcs used in acceptable paths must belong to either $M(a)_d$ or $M(a)_d^T$. By applying $\sigma_T^{-1}$ to the columns of M(s), we obtain a matrix whose diagonal elements all have the value zero, while all other entries have a positive value. This is because $\sigma_T^{-1}(M(s)) = \sigma_T^{-1}(M(s))^-$. It follows that all acceptable paths contain only positive values, implying that all acceptable cycles have positive values.